\documentclass[12pt]{article}

\usepackage[margin=1in]{geometry}  
\usepackage{graphicx}              
\usepackage{amsmath,euscript}               
\usepackage{amsfonts}              
\usepackage{amssymb,epsf,amsthm, epsfig}                
\usepackage[toc,page]{appendix}

\newtheorem{thm}{Theorem}[section]

\newtheorem{prop}[thm]{Proposition}

\newtheorem{lemma}[thm]{Lemma}

\theoremstyle{definition}
\newtheorem{rem}[thm]{Remark}
\newtheorem{defn}[thm]{Definition}

\input xy
\xyoption{all}

\newcommand{\Z}{\mathbb{Z}}

\newcommand{\calN}{\mathbb{N}}

\newcommand{\CC}{\EuScript C}
\newcommand{\DD}{\EuScript D}

\newcommand{\SSS}{\EuScript S}

\newcommand{\pf}{{\it Proof.}\hspace{2ex}}

\newcommand{\epf}{\hspace*{\fill}\mbox{$\halmos$}}
\newcommand{\halmos}{\rule{1ex}{1.4ex}}

\newcommand{\GL}{\mathop{\mathrm{GL}}\nolimits}

\newcommand{\Hom}{\mathop{\mathrm{Hom}}\nolimits}

\newcommand{\modu}{\mathop{\mathrm{mod}}\nolimits}

\newcommand{\sg}{\mathop{\mathrm{sg}}\nolimits}
\newcommand{\st}{\mathop{\mathrm{st}}\nolimits}
\newcommand{\qst}{\mathop{\mathrm{qst}}\nolimits}

\newcommand{\Mat}{\mathop{\mathrm{Mat}}\nolimits}

\newcommand{\sdeg}{\mathop{\mathrm{\leq_{\deg}}}\nolimits}
\newcommand{\trideg}{\mathop{\mathrm{\leq_{\bigtriangleup}}}\nolimits}

\begin{document}


\title{ Triangle Order $\leq_{\bigtriangleup}$ in Singular Categories}
\date{}
\author{Zhengfang WANG\thanks{zhengfang.wang@imj-prg.fr, Universit\'e Paris Diderot-Paris 7, Institut
de Math\'ematiques de Jussieu-Paris Rive Gauche CNRS UMR 7586, B\^atiment Sophie Germain, Case 7012,
75205 Paris Cedex 13, France} 
}

\maketitle

\begin{abstract}
Degeneration of modules is defined geometrically. Riedtmann and Zwara show that this degeneration is equivalent to the existence of a certain short exact sequence. Then Yoshino and independently Jensen, Su and Zimmermann generalised this notion to triangulated categories. We write $X\leq_\Delta Y$ if $X$ degenerates to $Y$. In this paper, we prove that $\leq_{\bigtriangleup}$ applied to the singular category $\DD_{\sg}(A)$ of a finite-dimensional $k$-algebra $A$ induces a partial order on the set of isomorphism classes of objects in $\DD_{\sg}(A)$.
\end{abstract}

\section{Introduction}
Degeneration order of modules is introduced from
geometric methods of representation theory of finite
dimensional algebras. More precisely, let $A$ be a finite dimensional associative $k$-algebra over the algebraically closed field $k$. Let $d$ be an positive integer. A $d$-dimensional (left) $A$-module $M$ is the vector space $k^d$ together with an action by $A$ from the left. We denote by $\modu_d(A)$ the set of $d$ dimensional $A$-modules. Note that $\modu_d(A)$ is an affine variety (For more ample details we refer to Section 2). The general linear group $\GL_d(k)$ acts on $\modu_d(A)$ by conjugation. The orbits under this action are the isomorphism classes of $d$-dimensional $A$-modules. We say that an $A$-module $N$ is called a degeneration of $M$ (denote by $M\leq_{\deg} N$) if $N$ belongs to the Zariski closure of the $\GL_d(k)$-orbit of $M$ in $\modu_d(A)$. Clearly, this degeneration defines a partial order on the set of isomorphism classes of $A$-modules.  Riedtmann and Zwara gave an algebraic description of degeneration order in \cite{Rie} and \cite{Zwa2}.  They showed that
$M\sdeg N$ if and only if there is an $A$-module $Z$ and an
exact sequence
\begin{eqnarray*}
  0\rightarrow N\rightarrow M\oplus Z\rightarrow Z\rightarrow 0
\end{eqnarray*}
and equivalently there exist an $A$-module $Z'$ and an exact sequence
\begin{eqnarray*}
  0\rightarrow Z'\rightarrow Z'\oplus M\rightarrow N\rightarrow 0.
\end{eqnarray*}
Later in \cite{Yosh2} Yoshino gave a scheme-theoretical definition
of degenerations, so that it can be considered for modules of  a
Noetherian algebra.
In this paper, we consider the algebraic description of degeneration order as our definition of degeneration order of $A$-$\modu$. That is, let $A$ be a
finite dimensional $k$-algebra $A$ over any field $k$ (not necessary algebraically closed), we define an $A$-module $N$ is a degeneration of an $A$-module $M$
 (still denoted by $\leq_{\deg}$) if there exists an $A$-module $Z$ and an exact sequence
\begin{eqnarray*}
  0\rightarrow N\rightarrow M\oplus Z\rightarrow Z\rightarrow 0.
\end{eqnarray*}
Then from Theorem 2.2 in \cite{Zwa1} it follows that
this degeneration $\leq_{\deg}$ is a partial order on
the set of isomorphism classes of $A$-modules.

Degeneration theory for triangulated categories and
derived categories (cf. \cite{JSZ1}, \cite{JSZ2}, \cite{SaZi}) has been studied. In a triangulated category we say for two objects $X$ and $Y$ that $X\leq_{\Delta} Y$ if there is an object $Z$ and a distinguished triangle $$Z[-1]\rightarrow N\rightarrow Z\oplus M\rightarrow Z.$$
In \cite{JSZ1}, Jensen, Su and Zimmermann showed that the triangle relation $\leq_{\bigtriangleup}$
in bounded derived category $\DD^b(A)$ is a partial order.  More generally, in \cite{JSZ2}, the authors
showed
that, under some finiteness assumptions on the triangulated category including the condition that the morphism spaces between objects are finite dimensional,
$\trideg$ is always a partial order.

The singular category was introduced by Ragnar-Olaf Buchweitz
in an unpublished manuscript \cite{Bu}, he called it stable
derived category, and Dmitri Orlov \cite{Orl} rediscovered
this notion independently in algebraic geometry and mathematical
physics, under the name of singular category. We remark that in general
the singular category $\DD_{\sg}(A)$  of a finite-dimensional
$k$-algebra $A$ is not Hom-finite and is in general not a
Krull-Schmidt category (cf. \cite{Chen}, \cite{ZhZi}).
Therefore we cannot use \cite{JSZ2} to argue
that the triangle relation $\trideg$ in $\DD_{\sg}(A)$ for
 any finite dimensional $k$-algebra $A$ is a partial order.
 But in this paper, we will prove, in a different way,  that $\trideg$ in $\DD_{\sg}(A)$  is really a partial order for
 any finite dimensional $k$-algebra $A$.

 This paper is organised as follows,
 in Section 2, we define a stable degeneration $\leq_{\st}$ for the stable
 module category $A$-$\underline{\modu}$ and prove that
 it is a partial order. In Section 3, we first recall some notions about the stabilization $\SSS(\mbox{$A$-$\underline{\modu}$})$ of the left triangulated category $A$-$\underline{\modu}$ and prove that in $\SSS(\mbox{$A$-$\underline{\modu}$})$, the triangle relation $\leq_{\bigtriangleup}$ coincides with the quasi-stable degeneration $\leq_{\qst}$ induced from the stable degeneration order $\leq_{\st}$ in $A$-$\underline{\modu}$.
We prove that the quasi-stable degeneration $\leq_{\qst}$ in $\SSS(\mbox{$A$-$\underline{\modu}$})$ is a partial order.  Last we show our main result (Theorem \ref{thm-main}), that is,  the triangle order $\trideg$
 in $\DD_{\sg}(A)$ is a partial order.
 We note that Theorem \ref{thm-main} and its proof extend without
 any changes to finitely generated modules over artinian algebras.

 For the background on the degeneration theory of modules and triangulated categories, we refer to
 \cite{JSZ1}, \cite{JSZ2}, \cite{Rie}, \cite{SaZi}, \cite{Yosh1}, \cite{Yosh2}, \cite{Zwa1} and \cite{Zwa2}. For the definition of singular category, we refer to \cite{Bu}, \cite{Orl} and \cite{Zim}.

\bigskip
\noindent
{\bf Acknowledgements:} This work is a part of author's PhD thesis, I would like to thank my PhD supervisor Alexander Zimmermann for his many valuable suggestions for improvement. I also want to thank Guodong Zhou for many inspiring conversations.

\section{Stable degeneration order}
First, let us recall the geometrical definition of degeneration order. Let $A$ be a finite dimensional associative $k$-algebra over an algebraically closed field $k$. Let $n$ be the dimension of $A$ and $d$ be an positive integer.  Let $\lambda_1, \lambda_2, \cdots, \lambda_n$ be a $k$-basis of $A$, $\lambda_i\lambda_j=\sum_l a_{ij}^l \lambda_l$ for $i, j=1, \cdots, n$ with the structure constants $a_{ij}^l\in k$. Then a module $M$ corresponds to a unique $n$-tuple of matrices $m=(m_1, \cdots, m_n)\in (\Mat_{d\times d}(k))^n$ such that $m_im_j=\sum_l a_{ij}^lm_l$ for $i, j=1, \cdots, n$. For each $1\leq l\leq n$ let $X^l$ denote the indeterminate matrix $(X_{\mu\nu}^l)_{\mu, \nu=1,\cdots, d}.$ Then there is a one-to-one correspondence between $\modu_d(A)$ and the zero set of the ideal $I\subset k[x_{\mu\nu}^l], (\mu, \nu=1,\cdots, d;  l=1, \cdots, n$), where $I$ is generated by the equations of the matrices $X^iX^j-\sum_la_{ij}^lX^l$ for $i, j=1,\cdots, n.$ The general linear group $\GL_d(k)$ acts on $\modu_d(A)$ by conjugation. The orbits under this action are the isomorphisms classes of $d$-dimensional $A$-modules. We say that an $A$-module $N$ is called a degeneration of $M$ (denote by $M\leq_{\deg} N$) if $N$ belongs to the Zariski closure of the $\GL_d(k)$-orbit of $M$ in $\modu_d(A)$. Clearly, this degeneration defines a partial order on the set of isomorphism classes of $A$-$\modu$. From the work of Riedtmann and Zwara in \cite{Rie} and \cite{Zwa2},  we have an algebraic description of degeneration order. In this paper, we use this algebraic description of degeneration order as our definition of degeneration order.
\begin{defn}\label{defn-deg}
  Let $k$ be a field. Let $A$ be a finite dimensional $k$-algebra and $X, Y\in \mbox{$A$-$\modu$}$,
  then we say that $X\leq_{\deg} Y$ if there exist $Z\in \mbox{$A$-$\modu$}$
  and an exact sequence on $A$-$\modu$,
  $$ 0\rightarrow Y\rightarrow X\oplus Z\rightarrow Z\rightarrow 0.$$
\end{defn}
\begin{rem}
From Theorem 2.3 in \cite{Zwa1}, $X\leq_{\deg}Y $ is equivalent to that there exists a short exact sequence
$$0\rightarrow Z'\rightarrow Z'\oplus X\rightarrow Y\rightarrow 0$$ for some $A$-module $Z'$.
  We remark that the degeneration order $\leq_{\deg}$ in $A$-$\modu$ defines a partial order on the set of isomorphism classes of $A$-modules (cf. \cite{Yosh2}, \cite{Zwa1}).
\end{rem}
\begin{defn}
  Let $A$ be a finite-dimensional $k$-algebra and $X, Y\in \mbox{$A$-$\modu$}$.
  We say that $X\leq_{\st}Y $ if and only if there exist two projective $A$-modules $P$ and
  $Q$ such that $X\oplus P\leq_{\deg} Y\oplus Q$. Clearly, we can induce a relation (called stable degeneration, still denote by $\leq_{\st}$) on isomorphism classes of objects
  in $A$-$\underline{\modu}$.
\end{defn}
\begin{rem}
Note that $P=0$ in $A$-$\underline{\text{mod}}$ if and only if $P$ is a projective module, so that $X\leq_{\st}Y$ is well-defined for any two objects $X$ and $Y$ in $A$-$\underline{\text{mod}}$.
\end{rem}

\begin{lemma}\label{lemma- 2.4}
  Let $A$ be a finite dimensional $k$-algebra. Then $\leq_{\st}$ defines a partial order on the set of isomorphism classes in $A$-$\underline{\modu}$.
\end{lemma}
\pf First, it's clear that the reflexivity is inherited from the degeneration order. We need to check anti-symmetry and transitivity.
For anti-symmetry, let $X\leq_{\st} Y$ and $Y\leq_{\st} X$. Therefore
there exist projective $A$-modules $P, Q, P', Q'$ such that
$$X\oplus P\leq_{\deg} Y\oplus Q$$
and $$Y\oplus P'\leq_{\deg} X\oplus Q'.$$
From the transitivity of $\leq_{\deg}$,
we have $$X\oplus P\oplus P'\leq_{\deg} X\oplus Q\oplus Q'.$$
Therefore we obtain (cf. Proposition 4.4 \cite{Yosh1})
\begin{eqnarray}\label{ineqna}
  \dim_k\Hom_A(P\oplus P', S)\leq \dim_k\Hom_A(Q\oplus Q', S)
\end{eqnarray}
for any simple $A$-module $S$. Note that there exists
a canonical bijection between the isomorphism classes of
projective indecomposable modules
 and the isomorphism classes of simple modules for a finite dimensional $k$-algebra $A$ (cf. e.g. \cite{Lein}), hence it follows
 from the inequality (\ref{ineqna}) that $P\oplus P'$ is a direct
summand of $Q\oplus Q'$. Since $X\oplus P\oplus P'\leq_{\deg} X\oplus Q\oplus Q',$
we obtain $\dim_k (P\oplus P')=\dim_k (Q\oplus Q')$, which implies $$P\oplus P'\cong Q\oplus Q'.$$
Since $X\oplus P\leq_{\deg} Y\oplus Q$ and $Y\oplus P'\leq_{\deg} X\oplus Q'$, we have
$$X\oplus P\oplus P'\leq_{\deg} Y\oplus Q\oplus P'$$
and $$Y\oplus P'\oplus Q\leq_{\deg} X\oplus Q'\oplus Q\cong X\oplus P\oplus P'.$$
From the anti-symmetry of the degeneration order $\leq_{\deg}$, we know
that $$X\oplus P\oplus P'\cong Y\oplus Q\oplus P'.$$ Hence we have $X\cong Y$ in
$A$-$\underline{\modu}$. This shows that $\leq_{\st}$ is anti-symmetric on the set of isomorphism classes of $A-\underline{\text{mod}}$.
In order to prove transitivity of $\leq_{\st}$, let $$X\leq_{\st} Y$$ and $$Y\leq_{\st} Z$$ in $A$-$\underline{\modu}$.
Then there exist projective $A$-modules $P, Q, R$ and $S$ such that
$$X\oplus P\leq_{\deg} Y\oplus Q$$ and
$$Y\oplus R\leq_{\deg} Z\oplus S.$$
Hence we have
$$X\oplus P\oplus R \leq_{\deg} Y\oplus Q\oplus R$$
and $$Y\oplus Q\oplus R\leq_{\deg} Z\oplus Q\oplus S.$$
From the transitivity of $\leq_{\deg}$ it follows that
$$X\oplus P\oplus R \leq_{\deg} Z\oplus Q\oplus S.$$
Hence $X\leq_{\st} Z$ in $A$-$\underline{\modu}$. Therefore the transitivity of $\leq_{\st}$ holds.
\epf

Next we define the triangle relation $\trideg$ for a (left) triangulated category. For the concept
of left triangulated categories we refer to \cite{BeMa} and \cite{KeVo}.
\begin{defn}[\cite{Yosh2}, \cite{JSZ2}]\label{defn-triangle-relation}
Let $\CC$ be a (left) triangulated category and
$X, Y\in \CC$, then we say that $X\leq_{\bigtriangleup} Y$ if there exist
$Z\in \CC$ and an exact triangle in $\CC$,
$$Z[-1]\rightarrow Y\rightarrow X\oplus Z\rightarrow Z.$$
\end{defn}
\begin{rem}
As well-known, $A$-$\underline{\modu}$ has a left triangulated structure with
the syzygy functor $\Omega_A$ as the translation functor (cf. \cite{BeMa}).
Hence, we can consider the triangle relation $\leq_{\bigtriangleup}$ in $A$-$\underline{\modu}$.
Next we will prove that $\leq_{\bigtriangleup}$ coincides with $\leq_{\st}$ in $A$-$\underline{\modu}$.
\end{rem}

\begin{lemma}\label{lemma-deg-equi}
Let $A$ be a finite dimensional $k$-algebra and let $X, Y\in \mbox{$A$-$\underline{\modu}$}$.
Then $X\leq_{\st} Y$ if and only if $X\leq_{\bigtriangleup} Y$ in $A$-$\underline{\modu}$.
\end{lemma}
\pf Suppose that $X\leq_{\st} Y$ in $A$-$\underline{\modu}$. Then
there exist projective $A$-modules $P$ and $Q$ such that $X\oplus P\leq_{\deg} Y\oplus Q$.
By Definition \ref{defn-deg}, there is an exact sequence in $A$-$\modu$,
$$0\rightarrow Y\oplus Q\rightarrow X\oplus P\oplus Z\rightarrow Z\rightarrow 0.$$
This induces the following exact triangle in $A$-$\underline{\modu}$,
$$\Omega_A(Z)\rightarrow Y\rightarrow X\oplus Z\rightarrow Z.$$
Hence $X\leq_{\bigtriangleup} Y$.

Conversely, assume that $X\leq_{\bigtriangleup} Y$. Then there exist an $A$-module $Z$
and an exact triangle,
$$\Omega_A(Z)\rightarrow Y\rightarrow X\oplus Z\rightarrow Z.$$
By the construction of the left triangulated structure in $A$-$\underline{\modu}$,
we know that there exist projective $A$-modules $P, Q$ and $R$ and an exact sequence in $A$-$\modu$,
$$0\rightarrow Y\oplus P\rightarrow X\oplus Z\oplus Q\rightarrow Z\oplus R\rightarrow 0.$$
Consider the following commutative diagram,
\begin{eqnarray*}
\xymatrix{
              &                             &                                       0                     &             0                &\\
            &                              &                            Z \ar@{=}[r]   \ar[u]          &                     Z  \ar[u]   &\\
0\ar[r] &   Y\oplus P\ar[r]    &    X\oplus Z\oplus Q\ar[r] \ar[u]^{\alpha}  &  Z\oplus R\ar[r] \ar[u]&  0\\
       0\ar[r]    &       Y\oplus P\ar[r]\ar@{=}[u]                        &                      \ker(\alpha)        \ar[u]\ar[r]              &              R\ar[u]  \ar[r] &0\\
       &  &  0\ar[u] &  0\ar[u]
}
\end{eqnarray*}
where $\alpha$ is the composition of $X\oplus Z\oplus Q\rightarrow Z\oplus R\rightarrow Z$.
Since $R$ is projective, the bottom row sequence splits, hence $\ker{\alpha}\cong Y\oplus P\oplus R$ and  we obtain the exact sequence,
\begin{eqnarray*}
  \xymatrix{
  0\ar[r] &  Y\oplus P\oplus R\ar[r] &   X\oplus Q\oplus Z\ar[r]^-{\alpha} &  Z\ar[r] & 0.
  }
\end{eqnarray*}
Therefore $X\oplus Q\leq_{\deg} Y\oplus P\oplus R$, hence $X\leq_{\st} Y$ in $A$-$\underline{\modu}$.
\epf

\begin{lemma}\label{lemma-deg-Omega}
Let $A$ be a finite dimensional $k$-algebra. If $X\leq_{\st} Y$,
then $\Omega_A(X)\leq_{\st} \Omega_A(Y)$.
\end{lemma}
\pf From Lemma \ref{lemma-deg-equi} it follows that $X\leq_{\st} Y$
if and only if $X\leq_{\bigtriangleup} Y$ in $A$-$\underline{\modu}$. Hence we have
the following exact triangle in $A$-$\underline{\modu}$,
$$\Omega_A(Z)\rightarrow Y\rightarrow X\oplus Z\rightarrow Z.$$
The Axioms of (left) triangulated categories imply that there is the following exact triangle,
$$\Omega^2_A(Z)\rightarrow \Omega_A(Y)\rightarrow \Omega_A(X)\oplus \Omega_A(Z)
\rightarrow \Omega_A(Z).$$
Therefore $\Omega_A(X)\leq_{\bigtriangleup} \Omega_A(Y)$, hence $\Omega_A(X)\leq_{\st}
\Omega_A(Y)$ using again Lemma \ref{lemma-deg-equi}.
\epf

\section{Stable categories and stabilization}
We recall some notions about the stabilization of stable categories. For details,
we refer to \cite{Bel}.
\begin{defn}
Let $(\CC, \Omega, \bigtriangleup)$ be a left triangulated category.
The stabilization of $\CC$ is a pair $(\iota, \SSS(\CC))$, where $\SSS(\CC)$ is
a triangulated category and $\iota: \CC\rightarrow \SSS(\CC)$ is an exact functor,
called the stabilization functor, such that for any exact functor $F:\CC\rightarrow \DD$
to a triangulated category $\DD$, there exists a unique exact functor $F^*:\SSS(\CC)
\rightarrow \DD$ such that $F^*\iota=F$.
\end{defn}

We recall the construction of $\SSS(\CC)$ (cf. \cite{Bel}, \cite{Hel}, \cite{KeVo}).
An object of $\SSS(\CC)$ is a pair $(X, m)$ where $X\in\CC$ and $m\in\Z$.
$$\Hom_{\SSS(\CC)}((X, m), (Y, n)):=\lim_{\substack{\longrightarrow\\k\geqslant m, n}}\Hom_{\CC}(\Omega^{k-m}X, \Omega^{k-n}Y).$$
Then $\SSS(\CC)$ is a triangulated category with the translation functor:
$\tilde{\Omega}:\SSS(\CC)\rightarrow\SSS(\CC)$ defined as follows:
$\tilde{\Omega}(X, m)=(X, m-1)$.
$$\tilde{\Omega}(Z, l)\rightarrow (X, m)\rightarrow (Y, n)\rightarrow (Z, l)$$
is an exact triangle in $\SSS(\CC)$ if and only if there exist $k\in 2\Z$
and a triangle, which represents the above triangle, $$\Omega(\Omega^{k-l}(Z))\rightarrow
\Omega^{k-m}(X)\rightarrow \Omega^{k-n}(Y)\rightarrow \Omega^{k-l}(Z)$$
in $\CC$.
\begin{thm}[Corollary 3.9. \cite{Bel}]\label{thm-bel}
Let $A$ be a finite-dimensional $k$-algebra,  then there exists a
triangle equivalence $$\SSS(\mbox{$A$-$\underline{\modu}$})\cong \DD_{\sg}(A).$$
\end{thm}

\begin{defn}\label{defn-stable-degeneration1}
Let $A$ be a finite dimensional $k$-algebra.
Let $(X, m), (Y, n)\in \SSS(\mbox{$A$-$\underline{\modu}$})$.  We say that  $(Y, n)$ is a quasi-stable degeneration of $(X, m)$ (denote by $(X, m)\leq_{\qst} (Y, n)$)
if and only if there exists $k\in \calN$ such that $\Omega^{k-m}(X)\leq_{\st} \Omega^{k-n}(Y)$ in $A$-$\underline{\modu}$.
\end{defn}
\begin{rem}
Since $\SSS(\mbox{$A$-$\underline{\modu}$})$ is a triangulated category,
there is the triangle relation $\leq_{\bigtriangleup}$ (cf. Definition \ref{defn-triangle-relation}) in $\SSS(\mbox{$A$-$\underline{\modu}$})$.
Next we will show that these two relations $\leq_{\qst}$ and $\leq_{\bigtriangleup}$ in $\SSS(\mbox{$A$-$\underline{\modu}$})$ coincide.
\end{rem}

\begin{prop}\label{prop-stab-order}
Let $A$ be a finite dimensional $k$-algebra, Then in $\SSS(\mbox{$A$-$\underline{\modu}$})$
 we have that $(X, m)\leq_{\qst} (Y, n)$ if and only if $(X, m)\leq_{\bigtriangleup} (Y, n)$.
\end{prop}
\pf If $(X, m)\leq_{\qst} (Y, n)$, then by Definition \ref{defn-stable-degeneration1},  there exists $k\in\calN$ such that $\Omega^{k-m}(X)
\leq_{\st} \Omega^{k-n}(Y)$ in $A$-$\underline{\modu}$, which means that there exist an $A$-module $Z$ and an exact triangle
in $A$-$\underline{\modu}$
$$\Omega(Z)\rightarrow \Omega^{k-n}(Y)\rightarrow \Omega^{k-m}(X)\oplus Z \rightarrow Z.$$
So we have a triangle in $\SSS(\mbox{$A$-$\underline{\modu}$})$
$$(Z, -1)\rightarrow (\Omega^{k-n}(Y), 0)\rightarrow (\Omega^{k-m}(X), 0)\oplus (Z, 0) \rightarrow (Z, 0),$$
which is isomorphic to
$$(Z, -1)\rightarrow (Y, n-k)\rightarrow (X, m-k)\oplus (Z, 0)\rightarrow (Z, 0).$$
So $(X, m-k)\leq_{\bigtriangleup} (Y, n-k),$ hence $(X, m)\leq_{\bigtriangleup} (Y, n).$

Conversely, suppose that $(X, m)\leq_{\bigtriangleup} (Y, n)$. Then there exist $(Z, r)\in \SSS(\mbox{$A$-$\underline{\modu}$})$
and an exact triangle,
$$(Z, r-1)\rightarrow (Y, n)\rightarrow (X, m)\oplus (Z, r)\rightarrow (Z, r),$$
that is, there exist $k\in 2\calN$ and an exact triangle in $A$-$\underline{\modu}$,
$$\Omega^{k-r+1}(Z)\rightarrow \Omega^{k-n}(Y)\rightarrow \Omega^{k-m}(X)\oplus \Omega^{k-r}(Z)
\rightarrow \Omega^{k-r}(Z).$$
Hence $\Omega^{k-m}(X)\leq_{\st}\Omega^{k-n}(Y)$, so $(X, m)\leq_{\qst} (Y, n)$.
\epf

Now let us prove the main theorem (cf. Theorem \ref{thm-main}). Before we come to the proof of Theorem \ref{thm-main}, we need the following lemma.

\begin{lemma}\label{lemma-last}
Let $A$ be a finite dimensional $k$-algebra. Then
the quasi-stable degeneration relation $\leq_{\qst}$ in $\SSS(\mbox{$A$-$\underline{\modu}$})$ is a partial order on the set of isomorphism classes of objects in
$\SSS(\mbox{$A$-$\underline{\modu}$})$.
\end{lemma}
\pf The reflexivity is inherited from the stable degeneration in $A$-$\underline{\modu}$.
For transitivity, let $(X, m)\leq_{\qst} (Y, n)$ and $(Y, n)\leq_{\qst} (Z, l)$, then
there exist  $k_1, k_2\in\calN$ such that
$$\Omega^{k_1-m}(X)\leq_{\st} \Omega^{k_1-n}(Y)$$
and $$\Omega^{k_2-n}(Y)\leq_{\st} \Omega^{k_2-l}(Z).$$
From Lemma \ref{lemma-deg-Omega}, we have that
$$\Omega^{k_1+k_2-m}(X)\leq_{\st} \Omega^{k_1+k_2-n}(Y)$$
and $$\Omega^{k_1+k_2-n}(Y)\leq_{\st} \Omega^{k_1+k_2-l}(Z).$$
By Lemma \ref{lemma- 2.4}, we know that the stable degeneration $\leq_{\st}$ in $A$-$\underline{\modu}$ is a partial order, hence
$$\Omega^{k_1+k_2-m}(X)\leq_{\st}\Omega^{k_1+k_2-l}(Z).$$
So $(X, m)\leq_{\qst} (Z, l)$, the transitivity holds in $\SSS(\mbox{$A$-$\underline{\modu}$}).$
Similarly, we can show that the anti-symmetric.
\epf

\begin{thm}\label{thm-main}
Let $A$ be a finite dimensional $k$-algebra, then the triangle order $\leq_{\bigtriangleup}$  is a partial order on the set of isomorphism classes of $\DD_{\sg}(A)$.
\end{thm}
\pf From Theorem \ref{thm-bel}, we know that $\DD_{\sg}(A)\cong \SSS(\mbox{$A$-$\underline{\modu}$})$ as triangulated categories.
Hence it's sufficient to show that the triangle order $\leq_{\bigtriangleup}$ in $\SSS(\mbox{$A$-$\underline{\modu}$})$
is a partial order. But this follows from Proposition \ref{prop-stab-order} and Lemma \ref{lemma-last}. Therefore
we have shown our result.
\epf


%

\bibliographystyle{plain}

\begin{thebibliography}{99}
%

%
%
\bibitem[Buch]{Bu}
Ragnar Olaf Buchweitz,
{\normalsize\it Maximal Cohen-Macaulay modules and Tate-cohomology over Gorenstein rings,}
manuscript Universit\"at Hannover 1986.
%
%
%
%
%
%
%
%
%
%
%
%
%
%
%
%
%
%
%
%
%
%
%
%
%
%


\bibitem[Bel]{Bel}
Apostolos Beligiannis,
{\normalsize\it The homological theory of contravariantly finite subcategories: Auslander-Buchweitz
contexts, Gorenstein categories and (co-) stabilization}
Communications in Algebra, {\bf 28} (10), 4547-4596 (2000).

\bibitem[BeMa]{BeMa}
Apostolos Beligiannis and Nikolaos Marmaridis,
{\normalsize\it Left triangulated categories arising from contravariantly finite subcategories,}
Communications in Algebra {\bf 22} (1994), no. 12, 5021-5036.


\bibitem[Chen]{Chen}
Xiao-Wu Chen,
{\normalsize\it The singular category of an algebra with
radical square zero,}
 Documenta Mathematical {\bf 16} (2011), no. 2-3, 299-212.

\bibitem[Hel]{Hel}
Alex Heller,
{\normalsize\it Stable homotopy categories,}
Bulletin of the American Mathematical Society {\bf 74} (1968), 28-63.


\bibitem[JSZ1]{JSZ1}
Bernt Tore Jensen, Xiuping Su and Alexander Zimmermann,
{\normalsize\it Degenerations for derived categories,}
Journal of Pure and Applied Algbra {\bf 198} (2005) 281-295.

\bibitem[JSZ2]{JSZ2}
Bernt Tore Jensen, Xiuping Su and Alexander Zimmermann,
{\normalsize\it Degeneration-like orders for triangulated
categories,} Journal of Algebra and its Applications {\bf 4} (2005) 587-597.




\bibitem[KeVo]{KeVo}
Bernhard Keller and Dieter Vossieck,
{\normalsize\it Sous les cat\'egories d\'eriv\'ees,}
Comptes Rendus de l'Acad\'emie des Sciences Paris, S\'erie I Math\'ematique
{\bf 305} (6) (1987) 225-228.


\bibitem[Lein]{Lein}
Tom Leinster,
{\normalsize\it The bijection between projective indecomposable and simple modules,}
arXiv:1410.3671 [math.RA].

\bibitem[Orl]{Orl}
Dmitri Orlov,
{\normalsize \it Derived categories of coherent
sheaves and triangulated categories of singularities,}
Algebra, arithmetic, and geometry: in honor of Yu. I. Manin.
Volume II, 503531, Progress in Mathematics, 270, Birkhuser
Boston, Inc., Boston, MA, 2009.

\bibitem[Rie]{Rie}
Christine Riedtmann,
{\normalsize\it Degenerations for representations of quivers with relations,}
Annales de l'\'Ecole Normale Sup\'erieure {\bf 19} (1986) 275-301.

\bibitem[SaZi]{SaZi}
Manuel Saorin and Alexander Zimmermann,
{\normalsize\it An axiomatic approach for degenerations in
triangulated categories,} Preprint (2014).



\bibitem[Yosh1]{Yosh1}
Yuji Yoshino,
{\normalsize\it On degeneration of Cohen-Macaulay modules,}
Journal of Algebra {\bf 248} (2002) 271-290.

\bibitem[Yosh2]{Yosh2}
Yuji Yoshino,
{\normalsize\it On degeneration of modules,}
Journal of Algebra {\bf 278} (2004) 217-226.

\bibitem[Yosh3]{Yosh3}
Yuji Yoshino,
{\normalsize\it Stable degeneration of Cohen-Macaulay modules,}
Journal of Algebra {\bf 332} (2011) 500-521.



\bibitem[ZhZi]{ZhZi}
Guo-Dong Zhou and Alexander Zimmermann,
{\normalsize\it On Singular Equivalence of Morita Type,}
Journal of Algebra {\bf 385} (2013) 64-79.


\bibitem[Zim]{Zim}
Alexander Zimmermann,
{\normalsize\it Representation Theory: A Homological Algebra Point of View,}
Springer Verlag London, 2014.



\bibitem[Zwa1]{Zwa1}
Grzegorz Zwara,
{\normalsize\it A degeneration-like order for modules,}
Archiv der Mathematik {\bf 71} (1998) 437-444.










\bibitem[Zwa2]{Zwa2}
Grzegorz Zwara,
{\normalsize\it Degenerations of finite dimension modules are given by extensions,}
Compositio Mathematica {\bf 121} (2000) 205-218.




\end{thebibliography}

\end{document}